\documentclass[11pt]{article}

\usepackage{amsmath,amsfonts,amssymb,amsthm,graphicx,graphics}

\usepackage{fancyhdr}
\newtheorem{theorem}{Theorem}[section]

\newtheorem{lem}[theorem]{Lemma}

 \newtheorem{deff}[theorem]{Definition}

\newtheorem{rem}[theorem]{Remark}

\newtheorem{cor}[theorem]{Corollary}

\newtheorem{thm}[theorem]{Theorem}

\def\provx#1#2#3#4{
\setbox1=\hbox{\kern1.5pt$\scriptstyle#3$}
\def\zeichen{#2}
\ifx\zeichen\empty\setbox0=\hbox to .75em{}\else\setbox0=\hbox
{\kern1.5pt$\scriptstyle#2$}\fi \dimen1=\dp0 \ifdim \dimen1=0pt
\advance \dimen1 by 1.5ex \else \advance \dimen1 by 1.2ex
\fi\dimen3=2ex\dimen4=.5ex\ifdim \wd0<\wd1 \dimen2=\wd1 \else
\dimen2=\wd0 \fi\hbox{$#1\hskip 5pt minus5pt\vrule height\dimen3
depth\dimen4\raise\dimen1\copy0\hskip-1\wd0 \lower\ht1
\copy1\hskip-1\wd1\vrule width\dimen2 height.7ex
depth-.6ex\hskip3pt minus1.5pt#4\hskip2pt plus2pt minus2pt$}}

\def\prov#1#2#3{
\setbox1=\hbox{\kern1.5pt$\scriptstyle#2$}
\def\zeichen{#1}
\ifx\zeichen\empty\setbox0=\hbox to .75em{}\else\setbox0=\hbox
{\kern1.5pt$\scriptstyle#1$}\fi \dimen1=\dp0 \ifdim \dimen1=0pt
\advance \dimen1 by 1.5ex \else \advance \dimen1 by 1.2ex
\fi\dimen3=2ex\dimen4=.5ex\ifdim \wd0<\wd1 \dimen2=\wd1 \else
\dimen2=\wd0 \fi\hbox{\hskip0pt plus 4pt $\vrule height\dimen3
depth\dimen4\raise\dimen1\copy0\hskip-1\wd0
\lower\ht1\copy1\hskip-1\wd1\vrule width\dimen2 height.7ex
depth-.6ex \hskip3pt minus1.5pt#3\hskip2pt plus2pt minus2pt$}}

\def\prv#1#2{
\setbox1=\hbox{\kern1.5pt$\scriptstyle#2$}
\ifx\zeichen\empty\setbox0=\hbox to .75em{}\else\setbox0=\hbox
{\kern1.5pt$\scriptstyle#1$}\fi \dimen1=\dp0 \ifdim \dimen1=0pt
\advance \dimen1 by 1.5ex \else \advance \dimen1 by 1.2ex
\fi\dimen3=2ex\dimen4=.5ex\ifdim \wd0<\wd1 \dimen2=\wd1 \else
\dimen2=\wd0 \fi\hbox{\hskip.5em$\vrule height\dimen3
depth\dimen4\raise\dimen1\copy0\hskip-1\wd0
\lower\ht1\copy1\hskip-1\wd1\vrule width\dimen2 height.7ex
depth-.6ex \hskip3pt minus1.5pt$}}

\mathchardef\str='1066
\def\negprov#1#2#3{
\setbox1=\hbox{\kern1.5pt$\scriptstyle#2$} \setbox4=\hbox{$\str$}
\def\zeichen{#1}
\ifx\zeichen\empty\setbox0=\hbox to 1em{}\else\setbox0=\hbox
{\kern1.5pt$\scriptstyle#1$}\fi \dimen1=\dp0 \ifdim \dimen1=0pt
\advance \dimen1 by 1.5ex \else \advance \dimen1 by 1.2ex
\fi\dimen3=2ex\dimen4=.5ex\ifdim \wd0<\wd1 \dimen2=\wd1 \else
\dimen2=\wd0
\fi\hbox{\hskip.5em$\kern-1.9pt\raise1pt\copy4\kern-\wd4\kern1.9pt\vrule
height\dimen3 depth\dimen4\raise\dimen1\copy0\hskip-1\wd0
\lower\ht1\copy1\hskip-1\wd1\vrule width\dimen2 height.7ex
depth-.6ex \hskip3pt minus1.5pt#3\hskip2pt plus2pt minus2pt$}}

\def\mod#1#2{
\def\zeichen{#1}
\hbox{\hskip 2pt plus3pt minus 2pt\vrule width.5pt height2ex
depth.5ex \vbox{\ifx\zeichen\empty\hbox to .75em{}\else
\hbox{\kern1.5pt $\scriptstyle#1$}\fi \kern2pt \hrule \kern1.7pt
\hrule\kern1.7pt} \hskip3pt minus 2pt$#2$}\hskip2pt plus3pt
minus2pt}

\def\notmod#1#2{\hbox{\hskip 2pt plus 3pt minus 3pt\vrule width.5pt
height2ex depth.5ex \vbox{\hbox{\kern1.5pt
$\scriptstyle#1$}\kern3pt \setbox0=\hbox{\kern2pt$\scriptstyle/$}
\hrule \kern-1.7pt \copy0 \kern-\ht0 \kern 1.7pt
\hrule\kern1.7pt}n \hskip3pt minus 2pt$#2$}\hskip2pt plus3pt
minus2pt}
\def\sq{\hbox{\rlap{$\sqcap$}$\sqcup$}}
\def\qed{\ifmmode\sq\else{\unskip\nobreak\hfil\penalty50\hskip1em\null
\nobreak\hfil\sq\parfillskip=0pt\finalhyphendemerits=0\endgraf}\fi\medskip}

\def\lleq{\hbox{\hskip3pt minus3pt\kern1pt\lower4pt
\vbox{\hbox{$\scriptstyle\ll$}
\kern-7pt\hbox{\kern1pt$\scriptstyle=$}}\hskip3pt minus 3pt}}

\mathchardef\res='1152 \mathchardef\qin='1062
\mathchardef\qprec='1036 \mathchardef\qless='474
\mathchardef\dpkt='72

\def\BPmessage{Proof Tree (bussproofs) style macros. Version 0.6c.}
\def\EnableBpAbbreviations{%
	\let\AX\Axiom 	\let\AXC\AxiomC 	\let\UI\UnaryInf
	\let\UIC\UnaryInfC 	\let\BI\BinaryInf 	 \let\BIC\BinaryInfC
	\let\TI\TrinaryInf 	\let\TIC\TrinaryInfC 	 \let\LL\LeftLabel
	\let\RL\RightLabel 	\let\DP\DisplayProof }

\def\ScoreOverhang{4pt}			
\def\ScoreOverhangLeft{\ScoreOverhang}
\def\ScoreOverhangRight{\ScoreOverhang}

\def\extraVskip{2pt}			
\def\ruleScoreFiller{\hrule}		

\def\defaultScoreFiller{\ruleScoreFiller}  
\def\defaultBuildScore{\buildSingleScore}  

\def\defaultHypSeparation{\hskip.2in}   

\def\labelSpacing{3pt}		

\def\proofSkipAmount{\vskip.8ex plus.8ex minus.4ex}

\def\theHypSeparation{\defaultHypSeparation}
\def\alwaysScoreFiller{\defaultScoreFiller}	
\def\alwaysBuildScore{\defaultBuildScore}
\def\theScoreFiller{\alwaysScoreFiller}	
\def\buildScore{\alwaysBuildScore}   
\def\hypKernAmt{0pt}	

\def\defaultLeftLabel{}
\def\defaultRightLabel{}

\def\myTrue{Y}
\def\bottomAlignFlag{N}
\def\centerAlignFlag{N}

\expandafter\ifx\csname newenvironment\endcsname\relax%
\def\makeatletter{\catcode`\@=11\relax}
\def\makeatother{\catcode`\@=12\relax}
\makeatletter
\def\newcount{\alloc@0\count\countdef\insc@unt}
\def\newdimen{\alloc@1\dimen\dimendef\insc@unt}
\def\newskip{\alloc@2\skip\skipdef\insc@unt}
\def\newbox{\alloc@4\box\chardef\insc@unt}
\makeatother \else
{\begin{center}\proofSkipAmount \leavevmode}%
{\DisplayProof \proofSkipAmount \end{center} } \fi

\def\thecur#1{\csname#1\number\theLevel\endcsname}

\newcount\theLevel    
\global\theLevel=0		
\newcount\myMaxLevel
\global\myMaxLevel=0
\newbox\myBoxA		
\newbox\myBoxB
\newbox\myBoxC
\newbox\myBoxD
\newbox\myBoxLL		
\newbox\myBoxRL
\newdimen\thisAboveSkip		
\newdimen\thisBelowSkip		
\newdimen\newScoreStart		
\newdimen\newScoreEnd
\newdimen\newCenter
\newdimen\displace
\newdimen\leftLowerAmt
\newdimen\rightLowerAmt
\newdimen\scoreHeight
\newdimen\scoreDepth

\setbox\myBoxLL=\hbox{\defaultLeftLabel}%
\setbox\myBoxRL=\hbox{\defaultRightLabel}%

\def\allocatemore{%
	\ifnum\theLevel>\myMaxLevel%
		\expandafter\newbox\curBox%
		\expandafter\newdimen\curScoreStart%
		\expandafter\newdimen\curCenter%
		\expandafter\newdimen\curScoreEnd%
		\global\advance\myMaxLevel by1%
	\fi%
}

\def\prepAxiom{%
	\advance\theLevel by1%
	\edef\curBox{\thecur{myBox}}%
	\edef\curScoreStart{\thecur{myScoreStart}}%
	\edef\curCenter{\thecur{myCenter}}%
	\edef\curScoreEnd{\thecur{myScoreEnd}}%
	\allocatemore%
}

\def\Axiom$#1\fCenter#2${%
	\prepAxiom%
	 \setbox\myBoxA=\hbox{$\mathord{#1}\fCenter\mathord{\relax}$}%
	\setbox\myBoxB=\hbox{$#2$}%
	\global\setbox\curBox=%
	     \hbox{\hskip\ScoreOverhangLeft\relax%
		 \unhcopy\myBoxA\unhcopy\myBoxB\hskip\ScoreOverhangRight\relax}%
	\global\curScoreStart=0pt \relax 	 \global\curScoreEnd=\wd\curBox
\relax 	\global\curCenter=\wd\myBoxA \relax
	\global\advance \curCenter by \ScoreOverhangLeft%
	\ignorespaces }

\def\AxiomC#1{		
	\prepAxiom%
	\setbox\myBoxA=\hbox{#1}%
	\global\setbox\curBox =%
		\hbox{\hskip\ScoreOverhangLeft\relax%
                        \unhcopy\myBoxA\hskip\ScoreOverhangRight\relax}%
        \global\curScoreStart=0pt \relax
        \global\curScoreEnd=\wd\curBox \relax
        \global\curCenter=.5\wd\curBox \relax
        \global\advance \curCenter by \ScoreOverhangLeft%
	\ignorespaces }

\def\prepUnary{%
	\ifnum \theLevel<1 		\errmessage{Hypotheses missing!}
	\fi%
	\edef\curBox{\thecur{myBox}}%
	\edef\curScoreStart{\thecur{myScoreStart}}%
	\edef\curCenter{\thecur{myCenter}}%
	\edef\curScoreEnd{\thecur{myScoreEnd}}%
}

\def\UnaryInf$#1\fCenter#2${%
	\prepUnary%
	\buildConclusion{#1}{#2}%
	\joinUnary%
	\resetInferenceDefaults%
	\ignorespaces%
}

\def\UnaryInfC#1{
	\prepUnary%
	\buildConclusionC{#1}%
	\joinUnary%
	\resetInferenceDefaults%
	\ignorespaces%
}

\def\prepBinary{%
	\ifnum\theLevel<2 		\errmessage{Hypotheses missing!}
	\fi%
	\edef\rcurBox{\thecur{myBox}}
	\edef\rcurScoreStart{\thecur{myScoreStart}}%
	\edef\rcurCenter{\thecur{myCenter}}%
	\edef\rcurScoreEnd{\thecur{myScoreEnd}}%
	\advance\theLevel by-1
	\edef\lcurBox{\thecur{myBox}}
	\edef\lcurScoreStart{\thecur{myScoreStart}}%
	\edef\lcurCenter{\thecur{myCenter}}%
	\edef\lcurScoreEnd{\thecur{myScoreEnd}}%
}

\def\BinaryInf$#1\fCenter#2${%
	\prepBinary%
	\buildConclusion{#1}{#2}%
	\joinBinary%
	\resetInferenceDefaults%
	\ignorespaces%
}

\def\BinaryInfC#1{%
	\prepBinary%
	\buildConclusionC{#1}%
	\joinBinary%
	\resetInferenceDefaults%
	\ignorespaces%
}

\def\prepTrinary{%
	\ifnum\theLevel<3 		\errmessage{Hypotheses missing!}
	\fi%
	\edef\rcurBox{\thecur{myBox}}
	\edef\rcurScoreStart{\thecur{myScoreStart}}%
	\edef\rcurCenter{\thecur{myCenter}}%
	\edef\rcurScoreEnd{\thecur{myScoreEnd}}%
	\advance\theLevel by-1
	\edef\ccurBox{\thecur{myBox}}
	\edef\ccurScoreStart{\thecur{myScoreStart}}%
	\edef\ccurCenter{\thecur{myCenter}}%
	\edef\ccurScoreEnd{\thecur{myScoreEnd}}%
	\advance\theLevel by-1
	\edef\lcurBox{\thecur{myBox}}
	\edef\lcurScoreStart{\thecur{myScoreStart}}%
	\edef\lcurCenter{\thecur{myCenter}}%
	\edef\lcurScoreEnd{\thecur{myScoreEnd}}%
}

\def\TrinaryInf$#1\fCenter#2${%
	\prepTrinary%
	\buildConclusion{#1}{#2}%
	\joinTrinary%
	\resetInferenceDefaults%
	\ignorespaces%
}

\def\TrinaryInfC#1{%
	\prepTrinary%
	\buildConclusionC{#1}%
	\joinTrinary%
	\resetInferenceDefaults%
	\ignorespaces%
}

\def\buildConclusion#1#2{
        \setbox\myBoxA=\hbox{$\mathord{#1}\fCenter\mathord{\relax}$}%
        \setbox\myBoxB=\hbox{$#2$}%
	\setbox\myBoxC =%
	      \hbox{\hskip\ScoreOverhangLeft\relax%
		 \unhcopy\myBoxA\unhcopy\myBoxB\hskip\ScoreOverhangRight\relax}%
	\newScoreStart=0pt \relax%
	\newCenter=\wd\myBoxA \relax%
	\advance \newCenter by \ScoreOverhangLeft%
	\newScoreEnd=\wd\myBoxC%
}

\def\buildConclusionC#1{
	\setbox\myBoxA=\hbox{#1}%
	\setbox\myBoxC =%
		\hbox{\hbox{\hskip\ScoreOverhangLeft\relax%
                        \unhcopy\myBoxA\hskip\ScoreOverhangRight\relax}}%
	\newScoreStart=0pt \relax%
	\newCenter=.5\wd\myBoxC \relax%
	\newScoreEnd=\wd\myBoxC%
        \advance \newCenter by \ScoreOverhangLeft%
}

\def\joinUnary{
	\global\advance\curCenter by -\hypKernAmt%
	\ifnum\curCenter<\newCenter%
		\displace=\newCenter%
		\advance \displace by -\curCenter%
		\kernUpperBox%
	\else%
		\displace=\curCenter%
		\advance \displace by -\newCenter%
		\kernLowerBox%
	\fi%
        \ifnum \newScoreStart < \curScoreStart %
		\global \curScoreStart = \newScoreStart \fi%
	\ifnum \curScoreEnd < \newScoreEnd %
		\global \curScoreEnd = \newScoreEnd \fi%
	\ifnum \curScoreStart<\wd\myBoxLL%
		\global\displace = \wd\myBoxLL%
		\global\advance\displace by -\curScoreStart%
		\kernUpperBox%
		\kernLowerBox%
	\fi%
	\buildScore%
	\buildScoreLabels%
	\global \setbox \curBox =%
		\vbox{\box\curBox%
			\vskip\thisAboveSkip \relax%
			\nointerlineskip\box\myBoxD%
			\vskip\thisBelowSkip \relax%
			\nointerlineskip\box\myBoxC}%
	\global \curScoreStart=\newScoreStart%
	\global \curScoreEnd=\newScoreEnd%
	\global \curCenter=\newCenter%
}
\def\kernUpperBox{%
		\global\setbox\curBox =%
			\hbox{\hskip\displace\box\curBox}%
		\global\advance \curScoreStart by \displace%
		\global\advance \curScoreEnd by \displace%
		\global\advance\curCenter by \displace%
}
\def\kernLowerBox{%
		\global\setbox\myBoxC =%
			\hbox{\hskip\displace\unhbox\myBoxC}%
		\global\advance \newScoreStart by \displace%
		\global\advance \newScoreEnd by \displace%
		\global\advance\newCenter by \displace%
}

\def\joinBinary{
	\setbox\myBoxA=\hbox{\theHypSeparation}
	\lcurScoreEnd=\rcurScoreEnd%
	\advance\lcurScoreEnd by\wd\lcurBox%
	\advance\lcurScoreEnd by\wd\myBoxA%
	\displace=\lcurScoreEnd%
	\advance\displace by -\lcurScoreStart%
	\lcurCenter=.5\displace%
	\advance\lcurCenter by\lcurScoreStart%
	\setbox\lcurBox=%
		\hbox{\box\lcurBox\unhcopy\myBoxA\box\rcurBox}%
	\displace=\newCenter%
	\advance\displace by -.5\newScoreStart%
	\advance\displace by -.5\newScoreEnd%
	\advance\lcurCenter by \displace%
	\edef\curBox{\lcurBox}%
	\edef\curScoreStart{\lcurScoreStart}%
	\edef\curScoreEnd{\lcurScoreEnd}%
	\edef\curCenter{\lcurCenter}%
	\joinUnary%
}

\def\joinTrinary{
	\setbox\myBoxA=\hbox{\theHypSeparation}
	\lcurScoreEnd=\rcurScoreEnd%
	\advance\lcurScoreEnd by\wd\lcurBox%
	\advance\lcurScoreEnd by\wd\ccurBox%
	\advance\lcurScoreEnd by2\wd\myBoxA%
	\displace=\lcurScoreEnd%
	\advance\displace by -\lcurScoreStart%
	\lcurCenter=.5\displace%
	\advance\lcurCenter by\lcurScoreStart%
	\setbox\lcurBox=%
		\hbox{\box\lcurBox\unhcopy\myBoxA\box\ccurBox%
				  \unhcopy\myBoxA\box\rcurBox}%
	\displace=\newCenter%
	\advance\displace by -.5\newScoreStart%
	\advance\displace by -.5\newScoreEnd%
	\advance\lcurCenter by \displace%
	\edef\curBox{\lcurBox}%
	\edef\curScoreStart{\lcurScoreStart}%
	\edef\curScoreEnd{\lcurScoreEnd}%
	\edef\curCenter{\lcurCenter}%
	\joinUnary%
}

\def\DisplayProof{%
	\ifnum \theLevel=1 \relax \else
		\errmessage{Proof tree badly specified.}%
	\fi%
	\edef\curBox{\thecur{myBox}}%
	\ifx\bottomAlignFlag\myTrue%
		\displace=0pt%
	\else%
		\displace=.5\ht\curBox%
		\ifx\centerAlignFlag\myTrue\relax
		\else%
		      	\advance\displace by -3pt%
		\fi%
	\fi%
	\leavevmode%
	\lower\displace\hbox{\copy\curBox}%
	\global\theLevel=0%
	\global\def\alwaysBuildScore{\defaultBuildScore}
	\global\def\alwaysScoreFiller{\defaultScoreFiller}
	\def\bottomAlignFlag{N} 	\def\centerAlignFlag{N}
	\resetInferenceDefaults%
	\ignorespaces }

\def\buildSingleScore{
	\displace=\curScoreEnd%
	\advance \displace by -\curScoreStart%
	\global\setbox \myBoxD =%
		\hbox to \displace{\expandafter\xleaders\theScoreFiller\hfill}%
}

\def\buildDoubleScore{
	\buildSingleScore%
	\global\setbox\myBoxD=%
		\hbox{\hbox to0pt{\copy\myBoxD\hss}\raise2pt\copy\myBoxD}%
}

\def\buildNoScore{
	\global\setbox\myBoxD=\hbox{\vbox{\vskip1pt}}%
}

\def\LeftLabel#1{%
	 \global\setbox\myBoxLL=\hbox{{#1}\hskip\labelSpacing}%
	\ignorespaces }
\def\RightLabel#1{%
	\global\setbox\myBoxRL=\hbox{\hskip\labelSpacing #1}%
	\ignorespaces }

\def\buildScoreLabels{%
	\scoreHeight = \ht\myBoxD%
	\scoreDepth = \dp\myBoxD%
	\leftLowerAmt=\ht\myBoxLL%
	\advance \leftLowerAmt by -\dp\myBoxLL%
	\advance \leftLowerAmt by -\scoreHeight%
	\advance \leftLowerAmt by \scoreDepth%
	\leftLowerAmt=.5\leftLowerAmt%
	\rightLowerAmt=\ht\myBoxRL%
	\advance \rightLowerAmt by -\dp\myBoxRL%
	\advance \rightLowerAmt by -\scoreHeight%
	\advance \rightLowerAmt by \scoreDepth%
	\rightLowerAmt=.5\rightLowerAmt%
	\displace = \curScoreStart%
	\advance\displace by -\wd\myBoxLL%
	\global\setbox\myBoxD =%
		\hbox{\hskip\displace%
			\lower\leftLowerAmt\copy\myBoxLL%
			\box\myBoxD%
			\lower\rightLowerAmt\copy\myBoxRL}%
	\global\thisAboveSkip = \ht\myBoxLL%
	\global\advance \thisAboveSkip by -\leftLowerAmt%
	\global\advance \thisAboveSkip by -\scoreHeight%
	\ifnum \thisAboveSkip<0 %
		\global\thisAboveSkip=0pt%
	\fi%
	\displace = \ht\myBoxRL%
	\advance \displace by -\rightLowerAmt%
	\advance \displace by -\scoreHeight%
	\ifnum \displace<0 %
		\displace=0pt%
	\fi%
	\ifnum \displace>\thisAboveSkip %
		\global\thisAboveSkip=\displace%
	\fi%
	\global\thisBelowSkip = \dp\myBoxLL%
	\global\advance\thisBelowSkip by \leftLowerAmt%
	\global\advance\thisBelowSkip by -\scoreDepth%
	\ifnum\thisBelowSkip<0 %
		\global\thisBelowSkip = 0pt%
	\fi%
	\displace = \dp\myBoxLL%
	\advance\displace by \rightLowerAmt%
	\advance\displace by -\scoreDepth%
	\ifnum\displace<0 %
		\displace = 0pt%
	\fi%
	\ifnum\displace>\thisBelowSkip%
		\global\thisBelowSkip = \displace%
	\fi
	\global\thisAboveSkip = -\thisAboveSkip%
	\global\thisBelowSkip = -\thisBelowSkip%
	\global\advance\thisAboveSkip by\extraVskip
	\global\advance\thisBelowSkip by\extraVskip
}

\def\resetInferenceDefaults{%
	\global\def\theHypSeparation{\defaultHypSeparation}%
	\global\setbox\myBoxLL=\hbox{\defaultLeftLabel}%
	\global\setbox\myBoxRL=\hbox{\defaultRightLabel}%
	\global\def\buildScore{\alwaysBuildScore}%
	\global\def\theScoreFiller{\alwaysScoreFiller}%
	\gdef\hypKernAmt{0pt}
}

\def\provx#1#2#3#4{
\setbox1=\hbox{\kern1.5pt$\scriptstyle#3$}
\def\zeichen{#2}
\ifx\zeichen\empty\setbox0=\hbox to .75em{}\else\setbox0=\hbox
{\kern1.5pt$\scriptstyle#2$}\fi \dimen1=\dp0 \ifdim \dimen1=0pt
\advance \dimen1 by 1.5ex \else \advance \dimen1 by 1.2ex
\fi\dimen3=2ex\dimen4=.5ex\ifdim \wd0<\wd1 \dimen2=\wd1 \else
\dimen2=\wd0 \fi\hbox{${#1}\hskip 5pt minus5pt\vrule height\dimen3
depth\dimen4\raise\dimen1\copy0\hskip-1\wd0 \lower\ht1
\copy1\hskip-1\wd1\vrule width\dimen2 height.7ex
depth-.6ex\hskip3pt minus1.5pt#4\hskip2pt plus2pt minus2pt$}}

\def\prov#1#2#3{
\setbox1=\hbox{\kern1.5pt$\scriptstyle#2$}
\def\zeichen{#1}
\ifx\zeichen\empty\setbox0=\hbox to .75em{}\else\setbox0=\hbox
{\kern1.5pt$\scriptstyle#1$}\fi \dimen1=\dp0 \ifdim \dimen1=0pt
\advance \dimen1 by 1.5ex \else \advance \dimen1 by 1.2ex
\fi\dimen3=2ex\dimen4=.5ex\ifdim \wd0<\wd1 \dimen2=\wd1 \else
\dimen2=\wd0 \fi\hbox{\hskip0pt plus 4pt $\vrule height\dimen3
depth\dimen4\raise\dimen1\copy0\hskip-1\wd0
\lower\ht1\copy1\hskip-1\wd1\vrule width\dimen2 height.7ex
depth-.6ex \hskip3pt minus1.5pt#3\hskip2pt plus2pt minus2pt$}}

\def\prv#1#2{
\setbox1=\hbox{\kern1.5pt$\scriptstyle#2$}
\ifx\zeichen\empty\setbox0=\hbox to .75em{}\else\setbox0=\hbox
{\kern1.5pt$\scriptstyle#1$}\fi \dimen1=\dp0 \ifdim \dimen1=0pt
\advance \dimen1 by 1.5ex \else \advance \dimen1 by 1.2ex
\fi\dimen3=2ex\dimen4=.5ex\ifdim \wd0<\wd1 \dimen2=\wd1 \else
\dimen2=\wd0 \fi\hbox{\hskip.5em$\vrule height\dimen3
depth\dimen4\raise\dimen1\copy0\hskip-1\wd0
\lower\ht1\copy1\hskip-1\wd1\vrule width\dimen2 height.7ex
depth-.6ex \hskip3pt minus1.5pt$}}

\mathchardef\str='1066
\def\negprov#1#2#3{
\setbox1=\hbox{\kern1.5pt$\scriptstyle#2$} \setbox4=\hbox{$\str$}
\def\zeichen{#1}
\ifx\zeichen\empty\setbox0=\hbox to 1em{}\else\setbox0=\hbox
{\kern1.5pt$\scriptstyle#1$}\fi \dimen1=\dp0 \ifdim \dimen1=0pt
\advance \dimen1 by 1.5ex \else \advance \dimen1 by 1.2ex
\fi\dimen3=2ex\dimen4=.5ex\ifdim \wd0<\wd1 \dimen2=\wd1 \else
\dimen2=\wd0
\fi\hbox{\hskip.5em$\kern-1.9pt\raise1pt\copy4\kern-\wd4\kern1.9pt\vrule
height\dimen3 depth\dimen4\raise\dimen1\copy0\hskip-1\wd0
\lower\ht1\copy1\hskip-1\wd1\vrule width\dimen2 height.7ex
depth-.6ex \hskip3pt minus1.5pt#3\hskip2pt plus2pt minus2pt$}}

\def\mod#1#2{
	\def\zeichen{#1} 	\hbox{\hskip 2pt plus3pt minus 2pt\vrule
width.5pt height2ex depth.5ex 		 \vbox{\ifx\zeichen\empty\hbox to
.75em{}\else 		\hbox{\kern1.5pt $\scriptstyle#1$}\fi 		 \kern2pt
		\hrule 		\kern1.7pt 		\hrule\kern1.7pt} 		 \hskip3pt minus
2pt$#2$}\hskip2pt 		plus3pt minus2pt}

\def\notmod#1#2{\hbox{\hskip 2pt plus 3pt minus 3pt\vrule width.5pt
		height2ex depth.5ex 		 \vbox{\hbox{\kern1.5pt
$\scriptstyle#1$}\kern3pt
		\setbox0=\hbox{\kern2pt$\scriptstyle/$} 		 \hrule 		 \kern-1.7pt
		\copy0 		\kern-\ht0 		\kern 1.7pt 		 \hrule\kern1.7pt}n
		\hskip3pt minus 2pt$#2$}\hskip2pt 		plus3pt minus2pt}
\def\sq{\hbox{\rlap{$\sqcap$}$\sqcup$}}
\def\qed{\ifmmode\sq\else{\unskip\nobreak\hfil\penalty50\hskip1em\null
\nobreak\hfil\sq\parfillskip=0pt\finalhyphendemerits=0\endgraf}\fi\medskip}

\def\lleq	{\hbox{\hskip3pt minus3pt\kern1pt\lower4pt
		\vbox{\hbox{$\scriptstyle\ll$}
		\kern-7pt\hbox{\kern1pt$\scriptstyle=$}}\hskip3pt minus 3pt}}

\mathchardef\res='1152 \mathchardef\qin='1062
\mathchardef\qprec='1036 \mathchardef\qless='474
\mathchardef\dpkt='72

\newcommand{\AC}{{\mathbf{AC}}}

\newcommand{\dbi}{\mbox{
${\mathbf\Delta}^{ 1}_{ 2}$--${\mathbf{CA}}+{\mathbf{BI}}$}}
\newcommand{\KP}{{\mathbf{KP}}}

\newcommand{\BI}{{\mathbf{BI}}}

\newcommand{\beq}{\begin{eqnarray}}
\newcommand{\eeq}{\end{eqnarray}}

\newcommand{\V}{V}
\newcommand{\RR}{R}

\newcommand{\LL}{{\mathbf{L}}}

\newcommand{\prf}{{\bf Proof\/}: }

\newcommand{\CH}{{\mathcal H}}

\newcommand{\und}{\,\wedge\,}

\newcommand{\psioi}{\psi_{\Omega}}

\newcommand{\TI}{{\mathrm{TI}}}

\newcommand{\ZFC}{{\mathbf{ZFC}}}
\newcommand{\ZF}{{\mathbf{ZF}}}
\newcommand{\CZF}{{\mathbf{CZF}}}

\newcommand{\KPP}{{\mathbf{KP}}({\mathcal P})}
\newcommand{\Deltaop}{\Delta_0^{\mathcal P}}

\newcommand{\Sigmap}{\Sigma^{\mathcal P}}
\newcommand{\Sigmaop}{\Sigma^{\mathcal P}}

\newcommand{\PC}{{\mathcal P}}

\newcommand{\RSOP}{RS_{\Omega}^{\mathcal P}}
\newcommand{\RSOPR}{RS_{\Omega}^{\mathcal P}(\mathsf{R})}

\newcommand{\POW}{{\mathbf{Pow}}}

\newcommand{\varthetai}{C}
\newcommand{\varphii}{A}
\newcommand{\MAC}{{\mathbf{MAC}}}
\newcommand{\ZZZ}{{\mathbf{Z}}}
\newcommand{\su}{\ell} 
\newcommand{\field}{{\mathrm{Field}}}
\newcommand{\Acck}{{\mathbf{Acc}}}
\newcommand{\Pow}{{\mathrm{Powerset}}}
\newcommand{\GAC}{\mathbf{AC}_{\!\mbox{\it\tiny global}}}

\begin{document}


\title{Power Kripke-Platek set theory and the axiom of choice}

\author{Michael Rathjen\\[1ex]
{\tiny Department of Pure Mathematics, University of Leeds}\\
{\tiny Leeds LS2 9JT, England, rathjen@maths.leeds.ac.uk}}
\date{}
\maketitle

\begin{abstract}{Whilst Power Kripke-Platek set theory, $\KPP$, shares many properties with ordinary Kripke-Platek set theory, $\KP$, in  several ways it behaves quite differently from $\KP$. This is perhaps most strikingly demonstrated by a result, due to Mathias, to the effect that adding the axiom of constructibility to $\KPP$ gives rise to a much stronger theory, whereas in the case of $\KP$ the constructible hierarchy provides an inner model, so that $\KP$ and
$\KP+V=L$ have the same strength.

This paper will be concerned with the relationship between $\KPP$ and  $\KPP$ plus the axiom of choice or even the global axiom of choice, $\GAC$.
Since $L$ is the standard vehicle to furnish a model in which this axiom holds, the usual argument for demonstrating that the addition of $\AC$ or $\GAC$ to $\KPP$ does not increase proof-theoretic strength does not apply in any obvious way. Among other tools, the paper uses techniques from ordinal analysis  to show that $\KPP+\GAC$ has the same strength as $\KPP$, thereby answering a question of Mathias. Moreover, it is shown that $\KPP+\GAC$
is conservative over $\KPP$ for $\Pi^1_4$ statements of analysis.

The method of ordinal analysis for theories with power set was developed in an earlier paper.
The technique allows one to compute witnessing information from infinitary proofs, providing
bounds for the  transfinite iterations of the power set operation that are provable in a theory.

As the theory $\KPP+\GAC$ provides a very useful tool for defining models and realizability models of other theories
 that are hard to construct without access to a uniform selection mechanism, it is desirable to determine its exact proof-theoretic strength. This knowledge can for instance be used to determine the strength of Feferman's operational set theory with power set operation as well as constructive Zermelo-Fraenkel set theory with the axiom of choice.
\\[2ex]
 \\
MSC Primary: 03F15 $\phantom{A}$ 03F05  $\phantom{A}$ 03F35  Secondary: 03F03
}\end{abstract}


\section{Introduction}
  A previous paper \cite{rathjen-KPP} gave a characterization
 of the smallest segment of the von Neumann hierarchy
 which is closed under the provable power-recursive functions of $\KPP$. It also furnished a proof-theoretic reduction of $\KPP$ to
 Zermelo set theory plus iterations of the powerset operation to any ordinal below the Bachmann-Howard ordinal.\footnote{The theories share the same $\Sigma_1^{\mathcal P}$ theorems, but are still distinct since
  Zermelo set theory does not prove $\Deltaop$-Collection whereas $\KPP$ does not prove full Separation.}
   The same bound also holds
 for the theory $\KPP+\AC$,
where $\AC$ stands for the {\em axiom of choice}.
These theorems
 considerably sharpen
 results
of  H. Friedman to the extent that $\KPP+\AC$ does
not prove the existence of the first non-recursive ordinal $\omega_1^{CK}$
(cf. \cite[Theorem 2.5]{friedmanpower} and \cite[Theorem 10]{mathias}).
However, \cite{rathjen-KPP} did not explicitly address the question of the proof-theoretic strength of $\KPP+\AC$.
In the present paper it will be shown that even adding global choice, $\GAC$, to $\KPP$ does not increase its proof-theoretic strength. This is in stark contrast to the axiom of constructibility $V=L$ which increases it as was shown by Mathias \cite{maclane}. That this increase is enormous will also been borne out by the results of this paper.

Since the global axiom of choice, $\GAC$, is less familiar, let us spell out the details. By $\KPP+\GAC$ we mean an extension of $\KPP$ where the language contains a new binary relation symbol $\RR$ and  the axiom schemes of
$\KPP$ are extended to this richer language and
the following axioms pertaining to $\RR$ are added:
\begin{eqnarray}\label{RR} (i) && \forall x\forall y\forall z[\RR(x,y)\wedge\RR(x,z)\to y=z]\\
    (ii) && \forall x[ x\ne \emptyset \to \exists y\in x\,\RR(x,y)].
    \end{eqnarray}

\section{Power Kripke-Platek set theory}
Before stating the axioms of $\KPP$, let us recall the axioms of $\KP$.
Roughly speaking,
 $\KP$ arises from $\ZF$ by completely omitting the power set axiom and
restricting separation and collection to set
bounded formulae but adding set induction (or
class foundation). Quantifiers of the forms $\forall x\in a$, $\exists x\in a$
are called {\em set
bounded}. {\em Set
bounded} or  {\em
$\Delta_0$-formulae} are
 formulae wherein all quantifiers are set
bounded. The axioms of $\KP$ consist of
{\em Extensionality, Pair, Union, Infinity, $\Delta_0$-Separation}
 $$\exists x\,\forall u\left[u\in x\leftrightarrow(u\in
 a\,\wedge\,\varphii(u))\right]$$
 for all $\Delta_0$-formulae
 $\varphii(u)$, {\em $\Delta_0$-Collection}
 $$\forall x\in a\,\exists y\,G(x,y)\,\to\,\exists z\,\forall x\in
 a\,\exists y\in z\,G(x,y)$$
 for all $\Delta_0$-formulae
$G(x,y)$, and
 {\em Set Induction}
 $$\forall x\,\left[(\forall y\in x\,\varthetai(y))\to
 \varthetai(x)\right]\,\to\,\forall x\,\varthetai(x)$$
 for all formulae $\varthetai(x)$.

A transitive set $A$ such that $(A,{\in})$ is a model of $\KP$ is
called an {\em admissible set}. Of particular interest are the
models of $\KP$ formed by segments of G\"odel's {\em constructible
hierarchy} $\LL$.  An ordinal $\alpha$ is {\em admissible} if the
structure $(\LL_{\alpha},{\in})$ is a model of $\KP$.

 $\KP$ is an
important set theory as a great deal of set theory requires only
the axioms of $\KP$ and its standard models,  the
admissible sets,  have been a major source of interaction between
model theory, recursion theory and set theory (cf. \cite{ba}).
{\em Power Kripke-Platek set theory} is obtained from $\KP$ by also viewing the creation of the powerset of any set
as a basic operation performed on sets. In the classical context, subsystems of $\ZF$ with Bounded Separation and Power Set have been studied by Thiele \cite{thiele}, Friedman \cite{friedmanpower}
and more recently in great depth by Mathias \cite{mathias}. They also  occur naturally in power recursion theory, investigated by Moschovakis \cite{mosch} and
Moss \cite{moss}, where one studies a notion of computability on the universe of sets which regards the
power set operation as an initial function.
Semi-intuitionistic set theories with Bounded    Separation
but containing
  the Power Set axiom were proposed by Pozsgay \cite{pozsgay1,pozsgay2} and then studied more
systematically by Tharp \cite{tharp}, Friedman \cite{frieda} and Wolf \cite{wolf}.
Such theories are naturally related to systems derived from topos-theoretic notions and to type theories (e.g., see \cite{ml}).
 Mac Lane has singled out and championed a particular fragment of  $\ZF$, especially in his book {\em Form and Function} \cite{maclane}. {\em Mac Lane Set Theory}, christened $\MAC$ in \cite{mathias}, comprises the axioms
 of Extensionality, Null Set, Pairing, Union, Infinity, Power Set, Bounded Separation, Foundation, and Choice.

To state the axioms of $\KPP$ it is convenient to introduce another type of bounded quantifiers.
\begin{deff}\label{deltap} {\em
We use subset bounded quantifiers $\exists x\subseteq y \;\ldots$ and $\forall x\subseteq y\;\ldots$ as abbreviations
for $\exists x(x\subseteq y\und \ldots)$ and $\forall x(x\subseteq y\to \ldots)$, respectively.

The $\Delta_0^{\mathcal P}$-formulae are the smallest class of formulae containing the atomic formulae closed under $\wedge,\vee,\to,\neg$ and the quantifiers
  $$\forall x\in a,\;\exists x\in a,\; \forall x\subseteq a,\; \exists x\subseteq a.$$
  A formula is in $\Sigma^{\mathcal P}$ if belongs to the smallest collection of formulae which contains the $\Delta_0^{\mathcal P}$-formulae and is closed under $\wedge,\vee$ and the quantifiers
  $\forall x\in a,\;\exists x\in a,\; \forall x\subseteq a$ and  $\exists x$.  A formula is $\Pi^{\mathcal P}$ if belongs to the smallest collection of formulae which contains the $\Delta_0^{\mathcal P}$-formulae and is closed under $\wedge,\vee$, the quantifiers
  $\forall x\in a,\;\exists x\in a,\; \forall x\subseteq a$ and  $\forall x$.
  }\end{deff}
\begin{deff}\label{KPP}{\em  $\KPP$ has the same language as $\ZF$. Its axioms are the following: Extensionality, Pairing, Union, Infinity, Powerset,
       $\Delta_0^{\mathcal P}$-Separation,    $\Delta_0^{\mathcal P}$-Collection and
       Set Induction (or Class Foundation).
       \\[1ex]
  The  transitive models of $\KPP$ have been termed {\bf power admissible} sets in \cite{friedmanpower}.
}\end{deff}

\begin{rem}{\em
 Alternatively, $\KPP$ can be obtained from $\KP$ by adding  a function symbol $\PC$ for the  powerset function as a primitive symbol to the language and the axiom
$$\forall y\,[y\in \PC(x)\leftrightarrow y\subseteq x]$$
and extending the schemes of $\Delta_0$ Separation and Collection to the $\Delta_0$-formulae of this new language.}
\end{rem}

\begin{lem}\label{weaker}
$\KPP$ is {\bf not} the same theory as $\KP+\POW$. Indeed, $\KP+\POW$ is a much weaker theory than $\KPP$
 in which one cannot prove the existence of $V_{\omega+\omega}$.
 \end{lem}
 \prf See \cite[Lemma 2.4]{rathjen-KPP}. \qed

\begin{rem}{\em The system $\KPP$ in the present paper is not quite the same as the theory $\KP^{\mathcal P}$ in Mathias'
paper \cite[6.10]{mathias}. The difference between $\KPP$ and $\KP^{\mathcal P}$  is that in the latter system
set induction only holds for $\Sigma_1^{\mathcal P}$-formulae, or what amounts to the same, $\Pi^{\mathcal P}_1$ foundation
($A\ne\emptyset \to \exists x\in A\;x\cap A=\emptyset$ for $\Pi_1^{\mathcal P}$ classes $A$).

Friedman \cite{friedmanpower} includes only Set Foundation in his formulation of a formal system $\mathbf{PAdm}^s$ appropriate to the concept of recursion in the power set operation $\mathcal P$.
}\end{rem}

\section{Extracting additional explicit results from ordinal analysis}
\cite{rathjen-KPP} featured an ordinal analysis of $\KPP$. As it turns out the technique can be augmented to also yield an ordinal analysis of $\KPP+\GAC$. The changes mainly concern the infinitary system $\RSOP$ of \cite[$\S5$]{rathjen-KPP}.
The modified infinitary  system $\RSOPR$ results from $\RSOP$ by the following changes:
\begin{itemize}
\item[(i)]  $\RSOPR$-terms are defined as in Definition 5.1 of \cite{rathjen-KPP}, except that in clause 3, $F(\vec x,y)$
is allowed to be any $\Deltaop$-formula of $\KPP+\GAC$, i.e., it may contain the relation symbol $\mathsf{R}$.
\item[(ii)] The axioms and rules of $\RSOP$ in  Definition 5.3 of \cite{rathjen-KPP}, have to be formulated with
respect to the richer language, i.e., $\Delta^{\mathcal P}_0$ and $\Sigma^{\mathcal P}$  refer to the language with $\mathsf{R}$ as a basic symbol. This affects (A1), (A3), (A6), (A7), and the rule $\Sigma^{\mathcal P}\mbox{-} Ref$. Moreover, to these axioms  one adds two new ones:
\begin{eqnarray} \mbox{(A8)} && \Gamma, \neg(\exists y\in t)(y\in t), (\exists y\in t)\,\mathsf{R}(t,y)
\\ \mbox{(A9)} && \Gamma, \neg \mathsf{R}(t,s), \neg \mathsf{R}(t,r), s=r\,.
\end{eqnarray}
\end{itemize}
With these changes, the embedding of $\KPP+\GAC$ into the infinitary proof system (cf. \cite[Theorem 6.9]{rathjen-KPP})
and  cut elimination in  $\RSOPR$  proceed in exactly the same way as for $\RSOP$ in
\cite[$\S7$]{rathjen-KPP}, yielding the following result:

\begin{cor}\label{Kore} Let $A$ be a $\Sigmaop$-sentence of $\KPP+\GAC$.
Suppose that $\KPP+\GAC\vdash A$. Then there exists an operator $\CH$ and an ordinal $\rho< \psioi(\varepsilon_{\Omega+1})$
  such that $$\provx{\CH}{\rho}{\rho}{A}\,.$$
  \end{cor}
  \prf The same proof as in \cite[Corollary 7.7]{rathjen-KPP} works here, too.
  $\CH$ and $\rho$ can be explicitly constructed from the proof of $A$ in $\KPP+\GAC$.
  \qed

  A refinement of \cite[Theorem 8.1]{rathjen-KPP} then yields partial conservativity of $\KPP+\GAC$ over $\KPP+\AC$.

  \begin{thm}\label{GACcons} Let $A$ be a
   $\Sigma^{\mathcal P}$ sentence of the language of set theory without $\mathsf{R}$.
   If $\KPP+\GAC\vdash A$ then  $\KPP+\AC\vdash A$.
   \end{thm}
  \prf Suppose $\KPP+\GAC\vdash \theta$, where $\theta$ is a
{ $\Sigmap$-sentence}.
 It follows from Corollary \ref{Kore} that one can explicitly find $\mathcal H$ and
$\tau<\psioi(\varepsilon_{\Omega+1})$ such that  $\provx{\CH}{\rho}{\rho}{A}\,.$
The refinement of \cite[Theorem 8.1]{rathjen-KPP} consists in noticing that if ${\mathcal H}(\emptyset) \subseteq C^{\Omega}(\omega_m(\omega^{\Omega+m}),0)$ and
$\Gamma(s_1,\ldots,s_n)$ is a sequent consisting only of $\Sigmaop$-formulae and
  $$\provx{\CH}{\alpha}{\rho}{\Gamma(s_1,\ldots,s_n)}$$
  with $\alpha,\rho<\Omega$,
  then, for all   variable assignments $\su:VAR\to V_{\psioi(\omega_m(\omega^{\Omega+m}))}$,
  $$ V_{\psioi(\omega_m(\omega^{\Omega+m}))}\models
  \Gamma(\su(s_1),\ldots,\su(s_n))\,,$$
  if the predicate $\mathsf{R}$ is interpreted as the graph of a choice function on  $V_{\psioi(\omega_m(\omega^{\Omega+m}))}$.
  The entire ordinal analysis can thus be carried out in $\KPP+\AC$ since this theory proves the existence of
  $V_{\psioi(\omega_m(\omega^{\Omega+m}))}$ as well as a choice function on this set (i.e. a function $f$ defined on this set satisfying $f(x)\in x$ whenever $x\ne\emptyset$).
  As $A$ is $\Sigma^{\mathcal P}$ and $\mathsf{R}$ does not occur in $A$, $A$ must be true. This shows that
  $\KPP+\AC\vdash A$.
  \qed

  The next result extracts further information from infinitary proofs in $\RSOP$.

\begin{thm}\label{3.1} If $\KPP+\AC\vdash \theta$, where $\theta$ is a
{ $\Sigmap$-sentence},
then one can explicitly find an ordinal (notation) $\tau<\psioi(\varepsilon_{\Omega+1})$
such that $$\KP+\AC+\mbox{\tt the von Neumann hierarchy $(V_{\alpha})_{\alpha\leq\tau}$ exists}\vdash \theta.$$
\end{thm}
\prf First note that $\AC$ can be formulated as a $\Pi^{\mathcal P}_1$-sentence.
From $\KPP+\AC\vdash \theta$ one obtains $\KPP\vdash \AC\to \theta$. As the latter statement
is equivalent to a { $\Sigmap$-sentence}, it follows from \cite[Corollary 7.7]{rathjen-KPP}
and a refinement of \cite[Theorem 8.1]{rathjen-KPP} (as above) that one can explicitly find
a $\tau<\psioi(\varepsilon_{\Omega+1})$ such that
\begin{eqnarray}\label{neu1} &&V_{\tau}\models\AC\to \theta.\end{eqnarray}
As before, the refinement of \cite[Theorem 8.1]{rathjen-KPP} consists in noticing that if ${\mathcal H}(\emptyset) \subseteq C^{\Omega}(\omega_m(\omega^{\Omega+m}),0)$ and
$\Gamma(s_1,\ldots,s_n)$ is a sequent consisting only of $\Sigmaop$-formulae and
  $$\provx{\CH}{\alpha}{\rho}{\Gamma(s_1,\ldots,s_n)}$$
  with $\alpha,\rho<\Omega$,
  then, for all   variable assignments $\su:VAR\to V_{\psioi(\omega_m(\omega^{\Omega+m}))}$,
  $$ V_{\psioi(\omega_m(\omega^{\Omega+m}))}\models
  \Gamma(\su(s_1),\ldots,\su(s_n))\,.$$
  The entire ordinal analysis can be carried out in $$\KP+\mbox{\tt the hierarchy $(V_{\delta})_{\delta<\tau}$ exists}$$ for a suitable $\tau$, e.g. $\tau =\psioi(\omega_m(\omega^{\Omega+m}))+\omega$.
  Thus $$\KP+\AC+\mbox{\tt the hierarchy $(V_{\delta})_{\delta<\tau}$ exists}\vdash\theta.$$
\qed
  Below we shall talk about well-orderings $\prec$.
  The field of $\prec$ is the set $\{u\mid\, \exists v\,(u\prec v\,\vee\,v\prec u)\}$.
  If $u\in\field(\prec)$ we denote by $\prec\restriction u$ the ordering $\prec$ restricted to the set
  $\{v\mid\,v\prec u\}$. $\prec\restriction  u$ is said to be an {\bf initial segment} of $\prec$.
  We say that $\prec$ is a {\bf cardinal} if there is no bijection between $\field(\prec)$ and
  $\{v\mid\,v\prec u\}$ for any $u\in\field(\prec)$.

\begin{thm}\label{3.2} Let $\tau$ be a limit ordinal. If $$\KP+\AC+\mbox{\tt the von Neumann hierarchy $(V_{\alpha})_{\alpha<\tau}$ exists}$$ proves a
$\Pi^1_4$ statements $\Phi$ of second order arithmetic, then
$$\ZZZ+\mbox{\tt the von Neumann hierarchy $(V_{\alpha})_{\alpha<\tau\cdot 4+4}$ exists}$$ proves $\Phi$.
\end{thm}
\prf We briefly recall the proof that $\ZF$ and $\ZFC$ prove the same $\Pi^1_4$ statements
of second order arithmetic. Assume $\ZFC\vdash \Phi$, where
$$\Phi\,=\,\forall x\subseteq\omega\exists y
\subseteq \omega \forall u\subseteq \omega\exists v\subseteq \omega\, \theta(x,y,u,v)$$
 with $\theta(x,y,u,v)$ an arithmetic formula. Now fix an arbitrary $x\subseteq \omega$ and
 build the relativized constructible hierarchy $L(x)$ which is a model of $\ZFC$, assuming $\ZF$ in the background. Then $L(x)\models  \forall u\subseteq \omega\exists v\subseteq \omega \theta(x,y,u,v)$ for some $y\subseteq \omega$ with $y\in L(x)$.
It then follows from a  version of Shoenfield's Absoluteness Lemma with subsets of $\omega$ as parameters
 that
$ \forall u\subseteq \omega\exists v\subseteq \omega \theta(x,y,u,v)$ holds in $V$ and hence $\Phi$
holds in $V$.

We would like to simulate the foregoing proof in the background theory  $$\ZZZ+\mbox{\tt the von Neumann hierarchy $(V_{\alpha})_{\alpha\leq\tau\cdot 4+4}$ exists}.$$
The idea is basically that for a given $x\subseteq \omega$ and a well-ordering $\prec$ one can simulate an initial segment of the  constructible hierarchy $L(x)$ along $\prec$.
Let's denote this by $L_{\prec}(x)$. $L_{\prec}(x)$ is basically a set of formal terms built
from the elements of the field of $\prec$. However, we also need to equip $L_{\prec}(x)$
 with an equivalence relation $\thickapprox$ such that $s\thickapprox s'$ signifies that $s$ and $s'$ denote the `same' set, and an elementhood relation $\varepsilon$ such that
 $s\varepsilon t \,\wedge\,s\thickapprox s'\to s'\varepsilon t$. This is a well known
 procedure in proof theory, so we shall not dwell on the details. For another approach based on the G\"odel functions see \cite{mathias}.

 We need a sufficiently long well-ordering to guarantee that $(L_{\prec}(x),\thickapprox,\varepsilon)$
 is a model of $T$. We resort to Hartog's construction (see e.g. \cite[7.34]{moschovakis}).
 Given a set $A$ let $<_{h(A)}$ be the well-ordering that arises from the set of all well-orderings on subsets
 of $A$ by singling out equivalence classes $$[\lessdot]_A=\{\vartriangleleft\,\mid\, \vartriangleleft\mbox{ is a well-ordering
  of a subset of $A$ order-isomorphic to $\lessdot$}\}$$
  and setting $$[\lessdot_1]_A<_{h(A)}[\lessdot_2]_A\,\mbox{ iff $\lessdot_1$ is isomorphic to a proper intial segment
  of $\lessdot_2$}.$$
  Since $\ZZZ$ has the powerset axiom and proves comparability of well-orderings,
  one can show that $<_{h(A)}$ is a well-ordering; moreover, $<_{h(A)}$ is the smallest well-ordering
  such that there is no injection of its field into $A$ (see \cite[7.34]{moschovakis}). As a result, there is no
  bijection between the field of $<_{h(A)}$ and the field of a proper initial segment of $<_{h(A)}$, that is,
  $<_{h(A)}$ behaves like a cardinal.

  Now, if $\lessdot\subseteq V_{\alpha}\times\V_{\alpha}$ then $\lessdot\in V_{\alpha+3}$ and hence $[\lessdot]_{V_{\alpha}}\in V_{\alpha+4}$. Therefore $<_{h(V_{\alpha})}$ is a well-ordering on
  a subset of $V_{\alpha+4}$ and hence $<_{h(V_{\alpha})}$ has a smaller order-type than
  $<_{h(V_{\alpha+4})}$, that is $<_{h(V_{\alpha})}$ is order-isomorphic to a proper initial segment of
  $<_{h(V_{\alpha+4})}$.
  Finally set \begin{eqnarray*} \prec &:=& <_{h(V_{\tau\cdot 4+\omega})}\\
  L^*(x) &:=&  (L_{\prec}(x),\thickapprox,\varepsilon).
  \end{eqnarray*}
  Our goal is to show that $L^*(x)$ is a model of
   $$T:=\KP+\AC+\mbox{\tt the von Neumann hierarchy $(V_{\alpha})_{\alpha<\tau}$ exists}.$$
  It is a straightforward matter to show that $L^*(x)$ satisfies Pairing, Union, Infinity, and $\Delta_0$-separation,
  similarly as one shows that any $L_{\lambda}$ with $\lambda$ a limit $>\omega$ satisfies these axioms.
 Next we address $\Delta_0$-Collection.  Note that $L^*(x)$ is a set with a definable  well-ordering. So we have definable Skolem functions on
  $L^*$. Moreover, for any $v\in\field(\prec)$ which is a limit one has a bijection between $L_{\prec\restriction v}(x)$ and $\prec\restriction v$. Now suppose $L^*(x)\models \forall u\in a\,\exists y\,\Psi(u,y,\vec b\,)$ with all parameters exhibited.
  Then there exists $v\in\field(\prec)$ of limit type such that $a,\vec b\in L_{\prec\restriction v}(x)$. Let $Y$ be the $\Sigma_1$ Skolem hull of $L_{\prec\restriction v}(x)$ in $L^*$. Note that there is a bijection between $Y$ and
   and $\field(\prec\restriction v)$. Then, similarly as in the  Condensation Lemma (e.g. \cite[II.5.2]{devlin})
  one proves that $Y$ is isomorphic to some $L_{\prec\restriction y}(x)$ with $v\preceq y$ and $y\in\field(\prec)$.
  The latter follows since $\prec$ is a cardinal. As a result, $L^*(x)\models \forall u\in a\,\exists y\in L_{\prec\restriction y}(x)\,\Psi(u,y,\vec b\,)$, showing $\Delta_0$-Collection.

  Finally we have to show that in $L^*$ the powerset operation can be iterated at least  $\tau$ times.
  By the above we know that all orderings $<_{h(V_{\alpha\cdot 4})}$ for $\alpha<\tau$ are isomorphic to initial segments  $\prec\restriction u_{\alpha}$ of $\prec$
  and that they  form an increasing sequence of cardinals.
  Using again a condensation argument (e.g.  \cite[II.5.5]{devlin}) one shows
  that $$L^*(x)\models {\mathcal P}(L_{\prec\restriction u_{\alpha}}(x))\subseteq  L_{\prec\restriction u_{\alpha+1}}(x)$$
  from which it follows that $$L^*\models \mbox{the von Neuman hierarchy $(V_{\beta})_{\beta<\tau}$ exists}.$$
  Since $L^*(x)$ also has a $\Delta_1$ definable well-ordering $\AC$ holds in $L^*(x)$. \qed

\begin{cor} If $\Phi$ is $\Pi^1_4$ sentence such that
$\KPP+\AC\vdash \Phi$ then $\KPP\vdash \Phi$.
\end{cor}
\prf Suppose $\KPP+\AC\vdash \Phi$ with $\Phi$ being $\Pi^1_4$.
Then $\Phi$ is also a  $\Sigmap$-sentence. By Theorem \ref{3.1}
one can explicitly find an ordinal (notation) $\tau<\psioi(\varepsilon_{\Omega+1})$
such that $$\KP+\AC+\mbox{\tt the von Neumann hierarchy $(V_{\alpha})_{\alpha\leq\tau}$ exists}\vdash \theta.$$
With Theorem \ref{3.2} we have
$$\ZZZ+\mbox{\tt the von Neumann hierarchy $(V_{\alpha})_{\alpha<\tau\cdot 4+4}$ exists}\vdash\Phi.$$
Now, $\KPP$ proves the existence of the ordinal $\tau\cdot 4+\omega$ and also
proves that $V_{\tau\cdot 4+\omega}$ is a model of $\ZZZ$, it follows that
$V_{\tau\cdot 4+\omega}\models\Phi$, whence $\Phi$ holds.\qed

\section{On Mathias' theory $\KP^{\mathcal P}$ and a miscellany of related work and questions}
As explained in the introduction, the difference between $\KPP$ and $\KP^{\mathcal P}$ (see \cite[6.10]{mathias})  is that in the latter system
set induction only holds for $\Sigma_1^{\mathcal P}$-formulae.

Let $\varphi$ be the usual two place Veblen function (see \cite{veb}) used by proof theorists.
It is by now well-known how one deals with  $\Sigma_1^{\mathcal P}$ set induction
proof-theoretically, employing partial cut elimination and an asymmetric interpretation of
quantifiers. Without the powerset function this is dealt with in \cite[\S3]{jaeger-habil}. With the powerset function it is carried out in \cite{cantini} and also follows from
\cite[6.1]{rathjen-prim} if one substitutes for $G$ the powerset function.
The upshot of this work is the following Theorem.

\begin{thm}\label{4.1}
If $\KP^{\mathcal P}+\AC\vdash \theta$, where $\theta$ is a
{ $\Sigmap$-sentence},
then one can explicitly find an ordinal (notation) $\tau<\varphi_{\omega}(0)$
such that $$\KP+\AC+\mbox{\tt the von Neumann hierarchy $(V_{\alpha})_{\alpha\leq\tau}$ exists}\vdash \theta.$$
\end{thm}

As a consequence, one arrives at the following:

\begin{cor} If $\Phi$ is $\Pi^1_4$ sentence such that
$\KP^{\mathcal P}+\AC\vdash \Phi$ then $\KP^{\mathcal P}\vdash \Phi$.
\end{cor}
\prf Suppose $\KP^{\mathcal P}+\AC\vdash \Phi$ with $\Phi$ being $\Pi^1_4$.
Then $\Phi$ is also a  $\Sigmap$-sentence. By Theorem \ref{4.1}
one can explicitly find an ordinal (notation) $\tau<\varphi_{\omega}(0)$
such that $$\KP^{\mathcal P}+\AC+\mbox{\tt the von Neumann hierarchy $(V_{\alpha})_{\alpha\leq\tau}$ exists}\vdash \theta.$$
With Theorem \ref{3.2} we have
$$\ZZZ+\mbox{\tt the von Neumann hierarchy $(V_{\alpha})_{\alpha<\tau\cdot 4+4}$ exists}\vdash\Phi.$$
Now, since $\KP^{\mathcal P}$ proves the existence of the ordinal $\tau\cdot 4+\omega$ and also
proves that $V_{\tau\cdot 4+\omega}$ is a model of $\ZZZ$, it follows that
$V_{\tau\cdot 4+\omega}\models\Phi$, whence $\Phi$ holds.
\qed

One can slightly strengthen the previous result by also adding $\Pi_1^{\mathcal P}$ set induction which is the same as $\Sigma_1^{\mathcal P}$ foundation.

Let $\Delta_0^{\mathcal P}\mbox{-}\mathrm{DC}$ be the schema saying that
whenever $\forall x\exists y R(x,y)$ holds for a $\Delta^{\mathcal P}$ predicate $R$ then for every set $z$ there exists a function $f$ with domain $\omega$ such that $f(0)=z$ and
$$ \forall n\in\omega \,R(f(n),f(n+1)).$$

\begin{lem}\label{DC} $\KP^{\mathcal P}+\AC\vdash \Delta_0^{\mathcal P}\mbox{-}\mathrm{DC}$.
\end{lem}
\prf We argue in $\KP^{\mathcal P}+\AC$. Assume  $\forall x\exists y R(x,y)$. In
$\KP^{\mathcal P}$ we have the function $\alpha\mapsto V_{\alpha}$. Fix a set $z$.
Using $\Delta_0^{\mathcal P}$ Collection we  get
$$\forall \alpha\exists \beta\,[z\in V_{\beta}\,\wedge\,\forall x\in V_{\alpha}\,\exists y\in V_{\beta}\,R(x,y)].$$
By taking the smallest $\beta$ we get a $\Sigma_1^{\mathcal P}$ class function $F$
such that
$$\forall \alpha[z\in V_{F(\alpha)}\,\wedge\,\forall x\in V_{\alpha}\,\exists y\in V_{F(\alpha)}\,R(x,y)].$$
Now let $\alpha_0:=0$, $\alpha_{n+1}=F(\alpha_n)$, and $\tau:=\sup_{n\in \omega}\alpha_n$.
Then $z\in V_{\tau}\,\wedge\,\forall x\in V_{\tau}\,\exists y\in V_{\tau}\,R(x,y)$.
 Using $\AC$ there is a choice function
$h:V_{\tau}\to V_{\tau}$ such that $R(x,h(x))$ holds for all $x\in V_{\tau}$.
Now define $f:\omega\to V_{\tau}$ by $f(0)=z$ and $f(n+1)=f(f(n))$.
Then $R(f(n),f(n+1))$ holds for all $n\in\omega$. \qed

\begin{lem} $\KP^{\mathcal P}+\AC\vdash \Pi_1^{\mathcal P}\mbox{ set induction}.$
\end{lem}
\prf We argue in $\KP^{\mathcal P}+\AC$. Suppose we have a counterexample to
$\Pi_1^{\mathcal P}\mbox{ set induction}.$
Then there is a $\Pi_1^{\mathcal P}$ predicate $Q(x)=\forall u\,S(u,x)$ with $S$ $\Delta_0^{\mathcal P}$ and a set $a_0$ such that
$\neg Q(a_0)$ and $\forall x\,[\forall y\in x\,Q(y)\to Q(x)]$.
The latter is equivalent to $\forall x\,[\neg Q(x) \to \exists y\in x\,\neg Q(y)]$ which yields
$$\forall x\,\forall u\,\exists y\,\exists v\,[\neg S(u,x) \to y\in x\,\wedge\,\neg S(y,v)].$$
The bracketed part of the latter statement is $\Delta_0^{\mathcal P}$.
 $\neg Q(a_0)$ entails that there exists $b_0$ such that $\neg S(a_0,b_0)$.
 Thus by Lemma \ref{DC} there exists a function $f$ with domain $\omega$ such that $f(0)=\langle a_0,b_0\rangle$
 and, for all $n\in\omega$, $f(n)$ is a pair $\langle a_n,b_n\rangle$ such that
 $$\neg S(a_n,b_n) \to a_{n+1}\in a_n\,\wedge\,\neg S(a_{n+1},b_{n+1}).$$
 By induction on $n$ we get $S(a_n,b_n) \,\wedge\,a_{n+1}\in a_n$ for all $n$, contradicting
 $\Delta_0^{\mathcal P}$ foundation.
 \qed

 \begin{cor} If $\Phi$ is $\Pi^1_4$ sentence such that
$\KP^{\mathcal P}+\AC+\Pi_1^{\mathcal P}\mbox{ set induction}\vdash \Phi$ then $\KP^{\mathcal P}\vdash \Phi$.
\end{cor}

Adrian Mathias \cite[Theorem 9]{mathias} proved that adding $V=L$ to $\KP^{\mathcal P}$
results in a stronger theory that proves the consistency of $\KP^{\mathcal P}$. This can be strengthened further.

\begin{thm}\label{stronger}\cite[Theorem 9]{mathias} $\KPP+V=L$ is much stronger than $\KPP$. Even
$\KP^{\mathcal P}+V=L$ is much stronger than $\KPP$.
\end{thm}
\prf We work in $\KP^{\mathcal P}+V=L$. Pick a limit ordinal $\kappa$ such that
${\mathcal P}(\omega)\in L_{\kappa}$. There exists a $\Sigma_1(L_{\kappa})$ map of $\kappa$ onto $L_{\kappa}$ (see \cite[II.6.8]{devlin}).
In particular there exists a map $f$ of $\kappa$ onto ${\mathcal P}(\omega)$.
We claim that $\kappa$ is uncountable, meaning that there is no surjective function that maps $\omega$ onto $\kappa$. This follows since otherwise there
would exist a surjection of $\omega$ onto the powerset of $\omega$, contradicting Cantor.

Next we aim at showing that $\KP^{\mathcal P}+V=L$ proves that
the Howard-Bachmann ordinal exists. To this end we invoke \cite{rathjen2005}
section 4.
Let $B$ be the set of ordinal notations and $<_{hb}$ be its ordering as defined in \cite[\S4]{rathjen2005}.
The class $\mathbf{Acc}$ of \cite[Definition 4.2]{rathjen2005} consists of those
ordinal notations  $a<_{hb}\Omega$ such that there exists
 an ordinal $\alpha$ and an order isomorphism $f$ between
  $\alpha$ and the initial segment of $B$ determined by $a$.
  However, since $B$ is a countable set such an $\alpha$ will always be
  $<\kappa$. As a result,
  $\Acck$ is actually a set in our background theory.
  One can then show that all the results in \cite[\S4]{rathjen2005} hold in
  our background theory.
   However, we can prove more.
  The class $\mathfrak M$ of \cite[Definition 4.7]{rathjen2005} will also be a set in our background theory and therefore
  the metainduction
  of \cite[Theorem 4.13]{rathjen2005} can be carried out as a formal induction in our background theory, and thus it can be
   show that the Bachmann-Howard ordinal exists as a set-theoretic ordinal.
   Consequently we can carry out the ordinal analysis of $\KPP$ from \cite{rathjen-KPP} inside our background theory.
  As a result, we get that $\KP^{\mathcal P}+V=L$ proves the
  $\Sigma^{\mathcal P}_1$ reflection principle for $\KPP$.
  \qed

  \begin{cor} $\KP^{\mathcal P}+V=L$ proves the
  $\Sigma^{\mathcal P}_1$ reflection principle for $\KPP$.
  \end{cor}

\section{Some applications}
The theory $\KPP+\GAC$ provides a very useful tool for defining models and realizability models of other theories
 that are hard to construct without access to a uniform selection mechanism.
 Since its exact proof-theoretic strength has been determined, this knowledge can be used to determine the strength of
 other theories, too. We give two examples, the first being
  Feferman's operational set theory with power set operation.

  J\"ager,  in \cite{jaeger-ost}, raised the following question at the very end of $\S4$ in connection with the strength of
  $\mathrm{OST}(\mathbb P)$.

  \begin{quote} Unfortunately, the combination of Theorem 18 and Theorem 22 does not
completely settle the question about the consistency strength of  $\mathrm{OST}(\mathbb P)$
yet. So far we have an interesting lower and an interesting upper bound,
but it still has to be determined what the relationship between $\KPP$ and
$\KPP + (V =L)$ is.
\end{quote}

The theory $\KPP+\AC$ has a translation into  $\mathrm{OST}(\mathbb P)$ by \cite[Theorem 18]{jaeger-ost} and by
\cite[Theorem 22]{jaeger-ost} $\mathrm{OST}(\mathbb P)$ can be interpreted in $\KPP+V=L$. However, by Theorem \ref{stronger} this leaves a huge gap.
The main difficulty is posed by the choice operator $\mathbb C$ of $\mathrm{OST}$.
The inductive definition of the relation ternary $R(a,b,c)$ that serves to interpret the application of $\mathrm{OST}(\mathbb P)$ in $\KPP$
is given in stages $R^{\alpha}$ in \cite[Definition 7]{jaeger-ost}. The assumption $V=L$ plays a central role in clause
22 when picking a $<_{L}$-least $c$ that fulfills the pertaining requirements. Here,  instead of $L_{\alpha}$ one can use $V_{\alpha}$, and instead of picking the $<_L$-least that works one can take the set $C$ of all possible candidates in $V_{\alpha}$ and then apply the global choice function to $C$ to select a particular one.
Thus together with \cite[Theorem 18]{jaeger-ost} we have
\begin{cor} The theories $\mathrm{OST}(\mathbb P)$ and $\KPP+\GAC$ are mutually interpretable in each other
and have the same strength as $\KPP$.\end{cor}
Similar results hold for versions of these three theories where the amount of $\in$-induction is restricted.
\\[2ex]
A question left open in \cite{acend} was that of the strength of constructive Zermelo-Fraenkel set theory with the axiom of choice. There $\CZF+\AC$ was interpreted in $\KPP+V=L$ (\cite[Theorem 3.5]{acend}). However, the (realizability) interpretation
works with $\GAC$ as well. Thus, as $\CZF+\AC$ proves  the power set axiom, it follows from \cite[Theorem 15.1]{ml} that
we have the following:

\begin{cor} $\CZF+\AC$ and $\KPP$ have the same strength.
\end{cor}

\section*{Acknowledgement} Thanks to Adrian Mathias for suggesting the problem.
The author would like to thank the Isaac Newton Institute for Mathematical Sciences, Cambridge, for support and hospitality during the programme `Mathematical, Foundational and Computational Aspects of the Higher Infinite' (HIF) in 2015, where work on this paper was undertaken.
Part of the material is based upon research supported by the EPSRC of the UK through grant No. EP/K023128/1. This research was also supported by a Leverhulme Research Fellowship.

\end{document}